\documentclass[11pt]{amsart}
\usepackage{mathrsfs}
\usepackage{amssymb}

\usepackage{titletoc}
\usepackage{amscd}
\usepackage{amsmath, amssymb}
\usepackage{amsfonts}
\usepackage[colorlinks,linkcolor=blue,citecolor=blue, pdfstartview=FitH]{hyperref}

\numberwithin{equation}{section}

\theoremstyle{plain}
\newtheorem{thm}{Theorem}[section]
\newtheorem{theorem}[thm]{Theorem}
\newtheorem{lemma}[thm]{Lemma}
\newtheorem{corollary}[thm]{Corollary}
\newtheorem{proposition}[thm]{Proposition}

\newtheorem{conjecture}[thm]{Conjecture}
\theoremstyle{definition}

\newtheorem{problem}[thm]{Problem}
\newtheorem{remark}[thm]{Remark}

\newtheorem{definition}[thm]{Definition}

\newtheorem{example}[thm]{Example}

\newtheorem{defn-thm}[thm]{Definition-Theorem}

\newcommand{\sO}{{\mathcal O}}

\newcommand{\C}{{\mathbb C}}

\newcommand{\N}{{\mathbb N}}
\renewcommand{\P}{{\mathbb P}}

\newcommand{\R}{{\mathbb R}}

\renewcommand{\S}{{\mathbb S}}

\newcommand{\Z}{{\mathbb Z}}

\newcommand{\qtq}[1]{\quad\mbox{#1}\quad}
\newcommand{\bp}{\bar{\partial}}
\newcommand{\Om}{\Omega}

\newcommand{\ts}{\otimes}

\newcommand{\btheorem}{\begin{theorem}}
\newcommand{\etheorem}{\end{theorem}}
\newcommand{\bproposition}{\begin{proposition}}
\newcommand{\eproposition}{\end{proposition}}
\newcommand{\bdefinition}{\begin{definition}}
\newcommand{\edefinition}{\end{definition}}
\newcommand{\bcorollary}{\begin{corollary}}
\newcommand{\ecorollary}{\end{corollary}}
\newcommand{\bproof}{\begin{proof}}
\newcommand{\eproof}{\end{proof}}
\newcommand{\bremark}{\begin{remark}}
\newcommand{\eremark}{\end{remark}}
\newcommand{\eexample}{\end{example}}
\newcommand{\bexample}{\begin{example}}
\newcommand{\la}{\langle}
\newcommand{\elemma}{\end{lemma}}
\newcommand{\blemma}{\begin{lemma}}
\newcommand{\ra}{\rangle}
\newcommand{\sq}{\sqrt{-1}}

\newcommand{\p}{\partial}

\renewcommand{\bar}{\overline}

\renewcommand{\phi}{\varphi}

\newcommand{\ee}{\end{eqnarray*}}
\newcommand{\be}{\begin{eqnarray*}}

\newcommand{\beq}{\begin{equation}}
\newcommand{\eeq}{\end{equation}}

\newcommand{\bd}{\begin{enumerate}}
\newcommand{\ed}{\end{enumerate}}

\renewcommand{\tilde}{\widetilde}

\renewcommand{\>}{\rightarrow}

\usepackage{fancyhdr}
\renewcommand{\leftmark}{{\scriptsize RC-positive metrics on rationally connected manifolds}}

\begin{document}
\title{RC-positive metrics on rationally connected manifolds} \makeatletter
\let\uppercasenonmath\@gobble
\let\MakeUppercase\relax
\let\scshape\relax
\makeatother
\author{Xiaokui Yang}
\date{}
\address{Morningside Center of Mathematics, Academy of Mathematics and\\ Systems Science, Chinese Academy of Sciences, Beijing, 100190, China}
\address{HCMS, CEMS, NCNIS, HLM, UCAS, Academy of Mathematics and\\ Systems Science, Chinese Academy of Sciences, Beijing 100190,
China} \email{\href{mailto:xkyang@amss.ac.cn}{{xkyang@amss.ac.cn}}}

\thanks{This work
was partially supported   by China's Recruitment
 Program of Global Experts and  NSFC 11688101.}

\maketitle
\begin{abstract} In this paper, we prove that if   a compact K\"ahler manifold  $X$
has a smooth Hermitian metric $\omega$ such that $(T_X,\omega)$ is
uniformly RC-positive, then $X$ is projective and rationally
connected. Conversely, we show that, if a projective manifold $X$ is
rationally connected, then the dual tautological line bundle
$\sO_{T_X^*}(-1)$ is uniformly RC-positive (which is equivalent to
the existence of some RC-positive complex Finlser metric on $X$). As
an application, we prove that if $(X,\omega)$ is a compact K\"ahler
manifold  with \emph{certain} quasi-positive holomorphic sectional
curvature, then $X$ is projective and rationally connected.
\end{abstract}



{\small{\setcounter{tocdepth}{1} \tableofcontents}}

\section{Introduction}

A projective manifold $X$ is called \emph{rationally connected} if
any two points of $X$ can be connected by some rational curves. It
is easy to show that on a rationally connected projective manifold,
one has
 $$ H^0(X,(T^*_X)^{\ts m})=0, \qtq{for every} m\geq 1.$$
A well-known conjecture of Mumford says that the converse is also
true.\\

\noindent\textbf{Conjecture.} Let $X$ be a projective manifold. If
$$  H^0(X,(T^*_X)^{\ts m})=0, \qtq{for all} m\geq 1,$$ then $X$ is
rationally connected.\\

\noindent This conjecture holds when $\dim X\leq 3$ (\cite{KMM92})
and not much has been known in higher dimensions, and we refer to
\cite{KMM92, Cam92, DPS96B, Kol96, GHS03, Pet06, BDPP13,CDP, CP14,
BC15, Cam16, LP17, CH17, Yang18} and the reference therein.\\

In this paper, we obtain a differential geometric criterion for
rational connectedness.

\bdefinition\label{Defintion} A Hermitian holomorphic vector bundle
$(\mathscr E,h)$ over a
 complex manifold $X$ is called \emph{uniformly RC-positive (resp. semi-positive)
 at  point $q\in X$}, if there exists some nonzero  vector $u\in T_qX$
such that for every nonzero vector $v\in \mathscr E_q$, \beq
R^\mathscr E(u,\bar u,v,\bar v)>0 \ \ \ (\text{resp.} \geq 0).\eeq
$(\mathscr E,h)$ is called \emph{uniformly RC-positive (resp.
semi-positive)} if it is uniformly RC-positive (resp. semi-positive)
 at each  point of $X$. $(\mathscr E,h)$ is called \emph{uniformly RC-quasi-positive}
 if it is uniformly RC-semi-positive at all  points of $X$, and uniformly RC-positive at some point of $X$.
\edefinition \noindent  The uniform RC-positivity in Definition
\ref{Defintion} is slightly stronger than the notion
``RC-positivity"
 introduced in \cite{Yang18}. For instance, we show that K\"ahler
 manifolds with positive holomorphic sectional curvature have
 uniformly RC-positive tangent bundle (see Theorem \ref{phscurc} and Remark \ref{remark}).

The first main result of this paper  is on deformations of
RC-quasi-positive Hermitian metrics, which demonstrates the
flexibility of this new concept.

\btheorem\label{quasi2} Let $X$ be a compact K\"ahler manifold.
Suppose there exists a smooth Hermitian metric $h$ on $T_X$ such
$(T_X,h)$ is uniformly RC-semi-positive over $X$. If there exist two
open
 subsets $S$ and $U$ of $X$ such that
\bd \item $U$ is strongly pseudoconvex and $\bar S\subset U$;

\item $(T_X,h)$ is uniformly RC-positive on $X\setminus S$.

\ed  Then  $X$ has a   uniformly RC-positive Hermitian metric
$\tilde h$. Moreover,  $X$ is a projective and rationally connected
manifold.\etheorem

\noindent The proof of Theorem \ref{quasi2} relies on some key
ingredients in our previous paper \cite{Yang18}, a conformal
perturbation method (Theorem \ref{perturb}) and the fundamental
result   of Campana-Demailly-Peternell  on the criterion of rational
connectedness. Recently, Lei Ni and Fangyang Zheng introduced  in
\cite{NZ1, NZ2} some interesting curvature notions which can ensure
the projectivity of compact K\"ahler manifolds. The following result
is a special case of Theorem \ref{quasi2}.

\bcorollary \label{main2} Let $X$ be a compact K\"ahler manifold. If
$X$ admits a smooth Hermitian metric $\omega$ such that
$(T_X,\omega)$ is uniformly RC-positive, then $X$ is a projective
and rationally connected manifold. \ecorollary

\bremark\label{finsler} Corollary \ref{main2} also holds if the
uniformly RC-positive Hermitian metric  is  generalized to a
uniformly RC-positive Finsler metric. See Theorem \ref{main20} for
more details. \eremark

\noindent As an application of Corollary \ref{main2}, we obtain the
following Liouville type result. \bcorollary \label{main02} Let $X$
be a compact K\"ahler manifold. If $X$ admits a smooth Hermitian
metric $\omega$ such that $(T_X,\omega)$ is uniformly RC-positive,
then there is  no non-constant holomorphic map from $X$ to a compact
complex manifold $Y$ where $Y$ is Kobayashi hyperbolic, or $Y$ has
nef cotangent bundle. \ecorollary

\noindent Corollary \ref{main02} may also hold if $X$ is
a compact complex manifold (see discussions in \cite{Yang18a,Yang18c}. \\

 It is well-known that a projective manifold $X$ is rationally
connected if and only if it has a very free rational curve $C$, i.e.
$T_X|_C$ is ample. Recall that, on a smooth curve $C$, a vector
bundle $\mathscr E$ is uniformly RC-positive if and only if
$\mathscr E$ is ample. Hence, roughly speaking, the uniform
RC-positivity of $T_X$ is an \emph{analytical analogue} of the
existence of very free rational curves on $X$ (see also
\cite[Theorem~6.6]{Pet12}). Before giving a converse to Corollary
\ref{main2}, we fix the notations. Let $\mathscr E$ be a holomorphic
vector bundle over $X$ and $\P(\mathscr E^*)$ be the projective
bundle of $\mathscr E$. The tautological line bundle of $\P(\mathscr
E^*)$ is denoted by $\sO_\mathscr E(1)$. For instance, $\mathscr E$
is called ample if $\sO_\mathscr E(1)$ is an ample line bundle. The
second main result of this paper is:

\btheorem\label{main3} Let $X$ be a projective manifold. If $X$ is
rationally connected, then the tautological dual line bundle
$\sO_{T_X^*}(-1)$ is uniformly
 RC-positive. \etheorem

\noindent  We also conjecture that a vector bundle $\mathscr E$ is
(uniformly) RC-positive if and only if $\sO_{\mathscr E^*}(-1)$ is
uniformly RC-positive (e.g. Conjecture \ref{conjecture}), which is
analogous to  a conjecture of Griffiths  that a vector bundle
$\mathscr E$ is Griffiths positive if and only if $\sO_\mathscr
E(1)$ is Griffiths positive. Hence, it is reasonable to expect that
rationally connected manifolds have (uniformly) RC-positive tangent
bundles. On the other hand, it is well-known that there exists a
one-to-one correspondence between the set of Hermitian metrics on
$\sO_{\mathscr E^*}(-1)$ and the set of Finsler metrics on $\mathscr
E^*$. Hence, by Corollary \ref{main2}, Remark \ref{finsler} and
Theorem \ref{main3}, we can deduce that
$X$ is rationally connected if and only if $X$ has certain ``RC-positive" Finsler metric.  \\

 As an application of Theorem \ref{quasi2}, we obtain

 \btheorem \label{QHSC} Let $(X,\omega)$ be a compact
K\"ahler manifold with  nonnegative holomorphic sectional curvature.
If there exist two open
 subsets $S$ and $U$ of $X$ such that $U$ is strongly pseudoconvex, $\bar S\subset U$ and $(X,\omega)$ has positive holomorphic sectional curvature on $X\setminus S$.
Then $X$ has a uniformly RC-positive Hermitian metric. In
particular, $X$ is a projective and rationally connected manifold.
 \etheorem

 \noindent In particular,

\bcorollary\label{PHSC} Let $(X,\omega)$ be a compact K\"ahler
manifold with positive holomorphic sectional curvature.  Then
$(T_X,\omega)$ is uniformly RC-positive. Moreover, $X$ is a
projective and rationally connected manifold. \ecorollary

\noindent Corollary \ref{PHSC} confirms a well-known conjecture
(\cite[Problem~47]{Yau82}) of  S.-T. Yau. It was firstly proved in
our previous paper \cite[Theorem~1.7]{Yang18}. The proof here is
slightly simpler than that in \cite{Yang18}. We also conjecture that

\begin{conjecture} Let $X$ be a compact K\"ahler manifold. If it has a smooth Hermitian metric $\omega$ with
quasi-positive holomorphic sectional curvature (or, with uniformly
RC-quasi-positive $(T_X,\omega)$), then $X$ is projective and
rationally connected.
\end{conjecture}

 The rest of the paper is organized as follows. In Section \ref{2},  we introduce the concept
 of uniform RC-positivity and investigate its geometric properties.
  In Section \ref{VRC}, we derive vanishing theorems for uniformly
RC-positive vector bundles and prove Theorem \ref{quasi2} and
Corollary \ref{main2}. In Section \ref{re}, we prove Theorem
\ref{main3} and propose more questions. In Section \ref{HSC00}, we
give a
 proof of Theorem \ref{QHSC}.  \\

\textbf{Acknowledgements.} I am very grateful to Professor
 K.-F. Liu and Professor S.-T. Yau for their support, encouragement and stimulating
discussions over  years. I would also like to thank Professors J.-Y.
Cao, J.-P. Demailly,  Huitao Feng, Y.-H. Wu, J. Xiao,
 and X.-Y. Zhou for helpful suggestions.

\vskip 2\baselineskip

\section{Uniformly RC-positive Hermitian vector bundles over complex manifolds }\label{2}

 Let $(\mathscr E,h)$ be a Hermitian holomorphic vector
bundle over a complex manifold $X$ with Chern connection $\nabla$.
Let $\{z^i\}_{i=1}^n$ be the  local holomorphic coordinates
  on $X$ and  $\{e_\alpha\}_{\alpha=1}^r$ be a local frame
 of $\mathscr E$. The curvature tensor $R^{\mathscr E}\in \Gamma(X,\Lambda^{1,1}T^*_X\ts \mathrm{End}(\mathscr E))$ has components \beq R^\mathscr E_{i\bar j\alpha\bar\beta}= -\frac{\p^2
h_{\alpha\bar \beta}}{\p z^i\p\bar z^j}+h^{\gamma\bar
\delta}\frac{\p h_{\alpha \bar \delta}}{\p z^i}\frac{\p
h_{\gamma\bar\beta}}{\p \bar z^j}.\label{cu}\eeq (Here and
henceforth we sometimes adopt the Einstein convention for
summation.) If $(X,\omega_g)$ is a  Hermitian manifold, then
$(T_X,g)$ has Chern curvature components \beq R_{i\bar j k\bar
\ell}=-\frac{\p^2g_{k\bar \ell}}{\p z^i\p\bar z^j}+g^{p\bar
q}\frac{\p g_{k\bar q}}{\p z^i}\frac{\p g_{p\bar \ell}}{\p\bar
z^j}.\eeq The Chern-Ricci curvature $\mathrm{Ric}(\omega_g)$ of
$(X,\omega_g)$ is represented by $R_{i\bar j}=g^{k\bar\ell} R_{i\bar
j k\bar\ell}$ and the second Chern-Ricci curvature
$\mathrm{Ric}^{(2)}(\omega_g)$ has components $R^{(2)}_{k\bar
\ell}=g^{i\bar j} R_{i\bar j k\bar\ell}.$

\bdefinition A Hermitian holomorphic vector bundle $(\mathscr E,h)$
over a complex manifold $X$ is called \emph{Griffiths positive} if
at each point $q\in X$ and for any nonzero vector $v\in \mathscr
E_q$, and any nonzero vector $u\in T_qX$,  \beq R^\mathscr E(u,\bar
u,v,\bar v)>0.\eeq  \edefinition

\noindent As analogous to Griffiths positivity, we introduced in
\cite{Yang18} the following concept.

 \bdefinition\label{Def}
A Hermitian holomorphic vector bundle $(\mathscr E,h)$ over a
complex manifold $X$ is called \emph{RC-positive
 at  point $q\in X$}, if for each nonzero vector $v\in \mathscr E_q$, there exists \textbf{some} nonzero  vector $u\in T_qX$
such that \beq R^\mathscr E(u,\bar u,v,\bar v)>0.\eeq  $(\mathscr
E,h)$ is called \emph{RC-positive} if it is  RC-positive
 at every  point of $X$.  \edefinition

\bremark Similarly, one can define semi-positivity, negativity and
etc.. For a Hermitian line bundle $(\mathscr L, h^{\mathscr L})$, it
is RC-positive if and only if its curvature
 $-\sq\p\bp\log h^{\mathscr L}$ has at least one positive
eigenvalue at each point of $X$.
\eremark

\noindent The following vanishing theorem is one of the key
ingredients by introducing the terminology ``RC-positivity".
\btheorem\label{vanishing} Let $X$ be a compact complex manifold. If
$(\mathscr E,h)$ is RC-positive, then \beq H^0(X,\mathscr
E^*)=0.\eeq \etheorem

\bproof An algebraic proof is included in
\cite[Lemma~2.10]{Yang18a}. Here we use an alternative proof by a
simple maximum principle. Since $(\mathscr E,h)$ is RC-positive, the
induced bundle $(\mathscr E^*,g)$ is RC-negative, i.e., at point
$q$, for any nonzero section $v$ of $\mathscr E^*$, there exists a
nonzero vector $u$ such that
$$R^{\mathscr E^*}(u,
\bar u, v,\bar v)<0$$
 For any $\sigma\in
H^0(X,\mathscr E^*)$, we have \beq \p\bp|s|^2_g=\la\nabla s, \nabla
s\ra_g-R^{\mathscr E^*}(\bullet, \bullet, s,\bar s).\label{max}\eeq
Suppose $|s|^2_g$ attains its maximum at some point $p$ and
$|s|^2_g(p)>0$. By applying maximum principle to (\ref{max}), we get
a contradiction. Hence, we deduce $s=0$ and $H^0(X,\mathscr E^*)=0$.
\eproof

\noindent In particular, we obtain a simple criterion for the
projectiveness of compact K\"ahler manifolds.

\bcorollary\label{kp} Let $X$ be a compact K\"ahler manifold.
Suppose $\Lambda^2T_X$ is $RC$-positive, then $X$ is projective.
\ecorollary

\bproof By Theorem \ref{vanishing}, we have
$H^{2,0}_{\bp}(X)=H^{0,2}_{\bp}(X)=0$. Hence, by the Kodaira theorem
(\cite[Theorem~1]{kod54}, see also \cite[Proposition~ 3.3.2 and
Corollary~5.3.3]{Huy05}), the K\"ahler manifold $X$ is projective.
\eproof

\noindent We introduce a notion slightly stronger than
RC-positivity.

 \bdefinition\label{Def1}
A Hermitian holomorphic vector bundle $(\mathscr E,h)$ over a
complex manifold $X$ is called \emph{uniformly RC-positive
 at  point $q\in X$}, if there exists some  vector $u\in T_qX$
such that for any nonzero vector $v\in \mathscr E_q$, one has  \beq
R^\mathscr E(u,\bar u,v,\bar v)>0.\eeq $(\mathscr E,h)$ is called
\emph{uniformly RC-positive} if it is uniformly RC-positive
 at every  $q\in X$.
 \edefinition

\bremark We can define uniform RC-negativity (resp. uniform
RC-nonnegativity, uniform RC-nonpositivity) in a similar way. We can
also define \emph{uniformly RC-positive along $k$-linearly
independent directions}, if there exist $k$ linearly independent
vectors $u_1,\cdots, u_k\in T_qX$ such that for any nonzero vector
$v\in \mathscr E_q$ and for each $i=1,\cdots, k$, one has $$
R^\mathscr E(u_i,\bar u_i,v,\bar v)>0.$$

  \eremark

\bremark It is easy to see that, for a line bundle $(\mathscr L,h)$,
RC-positivity and uniform RC-positivity are actually equivalent.

\eremark

\bproposition\label{est2} Let  $(\mathscr E,h)$ be a Hermitian
holomorphic vector bundle over a compact complex manifold $X$. Then
the following statements are equivalent:

\bd \item $(\mathscr E,h)$ is uniformly RC-positive;

\item for any Hermitian metric $\omega$ on $X$, there
exists a positive constant $C=C(\omega, h)$ such that: for any point
$q\in X$, there exists a unit vector $u\in T_qX$ such that \beq
R^\mathscr E(u,\bar u, v,\bar v)\geq C|v|_h^2, \qtq{for every} v\in
\mathscr E_q. \label{est} \eeq \ed

\bproof $(2)\Longrightarrow (1)$ is obvious. For $(2)\Longrightarrow
(1)$,  let \beq C=\inf_{q\in X} \sup_{u\in T_q X\setminus\{0\}}
\inf_{v\in \mathscr E_q\setminus\{0\}} \frac{R^\mathscr E(u,\bar
u,v,\bar v)}{|u|^2_\omega|v|_h^2}.\label{C}\eeq We claim $C>0$.
Indeed, if $C\leq 0$, by the compactness of $X$, there exists some
$q\in X$ such that
$$\sup_{u\in T_q X\setminus\{0\}} \inf_{v\in \mathscr E_q\setminus\{0\}} \frac{R^\mathscr E(u,\bar u,v,\bar
v)}{|u|^2_\omega|v|_h^2}\leq 0.$$  Hence, for any unit vector $u\in
T_q X$, $$\inf_{v\in \mathscr E_q\setminus\{0\}} \frac{R^\mathscr
E(u,\bar u,v,\bar v)}{|v|_h^2}\leq 0,$$ which contradicts to the
fact that $(\mathscr E,h)$ is uniformly RC-positive at  point $q\in
X$. It is obvious that (\ref{est}) follows from the definition
(\ref{C}) of the constant $C$. \eproof \eproposition

\bproposition\label{rcg} If $(\mathscr E,h_1)$ is uniformly
RC-positive and $(\mathscr F,h_2)$ is Griffiths semi-positive, then
$(\mathscr E\ts\mathscr F, h_1\ts h_2)$ is uniformly RC-positive.
\eproposition \bproof It follows from the curvature formula $
R^{\mathscr E\ts \mathscr F}=R^\mathscr E\ts \mathrm{Id}_\mathscr
F+\mathrm{Id}_\mathscr E\ts R^\mathscr F.$ 
 \eproof

\noindent The following  results can be deduced by similar methods
as in \cite[Theorem~3.5]{Yang18}.

 \bcorollary\label{relation} Let $(\mathscr E,h)$ be a Hermitian vector bundle over
a compact complex manifold $X$. \bd \item If $(\mathscr E,h)$ is
uniformly RC-positive, then $(\mathscr E,h)$ is RC-positive. \item
$(\mathscr E,h)$ is uniformly RC-positive if and only if $(\mathscr
E^*,h^*)$ is uniformly RC-negative;

 \item  If $(\mathscr E,h)$ is uniformly RC-negative, every subbundle $\mathscr S$ of $\mathscr E$ is uniformly RC-negative;

\item If $(\mathscr E,h)$ is uniformly RC-positive,
 every quotient bundle $\mathscr Q$ of $\mathscr E$ is uniformly RC-positive;

\item  If  $(\mathscr E,h)$ is uniformly RC-positive, every line subbundle $\mathscr L$ of
$\mathscr E^*$ is not pseudo-effective.

\ed\ecorollary




\noindent The uniform RC-positivity  has some important functorial
properties.

\bproposition\label{tensor} If $(\mathscr E,h)$ is  uniformly
RC-positive, then $(\mathscr E^{\ts m},h^{\ts m})$ are uniformly
RC-positive for $m\in\N^+$. Similarly, $\mathrm{Sym}^{\ts k}\mathscr
E$ $(k\in \N^+)$ and $\Lambda^p \mathscr E$ ($1\leq p\leq
\mathrm{rk}(\mathscr E)$) are all uniformly RC-positive.
\eproposition

\bproof  Fix a smooth metric on the compact complex manifold $X$.
For any $q\in X$, we choose a unit vector $u\in T_qX$ such that \beq
\inf_{v\in \mathscr E_q\setminus\{0\}}\frac{R^\mathscr E(u,\bar
u,v,\bar v)}{|v|_h^2}\geq C>0.\label{key}\eeq Let $\{e_1,\cdots,
e_r\}$ be a local unitary frame of $\mathscr E$ at point $q$ with
respect to $h$. Then for any local vector $v=\sum_{i_1,\cdots,
i_m}v^{i_1\cdots i_m}e_{i_1}\ts\cdots \ts e_{i_m}\in \Gamma(X,
\mathscr E^{\ts m}),$ by using the tensor product curvature formula
of $(\mathscr E^{\ts m},h^{\ts m})$, we have \be  R^{\mathscr E^{\ts
m}}(u,\bar u,v,\bar v) &=&\sum_{i_2,\cdots, i_m} R^{\mathscr
E}\left(u,\bar u,\sum_{i_1}v^{i_1i_2\cdots i_m}
e_{i_1},\bar{\sum_{j_1}v^{j_1i_2\cdots i_m}
e_{j_1}}\right)+\cdots\\&&+\sum_{i_1,i_2,\cdots, i_{m-1}}
R^{\mathscr E}\left(u,\bar u,\sum_{i_m}v^{i_1\cdots i_{m-1} i_m}
e_{i_m},\sum_{j_m}\bar{v^{i_1\cdots i_{m-1}j_m}
e_{j_m}}\right)\\
&\geq& C\left(\sum_{i_2,\cdots, i_m}\sum_{i_1}|v^{i_1i_2\cdots
i_m}|^2+\cdots+\sum_{i_1,\cdots,
i_{m-1}}\sum_{i_m}|v^{i_1i_2\cdots i_m}|^2\right)\\
& \geq& mC|v|^2,\ee where the first inequality follows from
(\ref{key}). Hence $(\mathscr E^{\ts m}, h^{\ts m})$ is uniformly
RC-positive. Similarly, we can show $\mathrm{Sym}^{\ts k}\mathscr E$
$(k\in \N^+)$ and $\Lambda^p \mathscr E$ ($1\leq p\leq
\mathrm{rk}(\mathscr E)$) are all uniformly RC-positive.\eproof

\bremark It is not hard to see that all Schur powers of a uniformly
RC-positive vector bundle are uniformly RC-positive. \eremark

\vskip 2\baselineskip

\section{Vanishing theorems and rational connectedness of compact K\"ahler
manifolds}\label{VRC}

In this section, we derive vanishing theorems for uniformly
RC-positive vector bundles and prove Theorem \ref{quasi2} and
Corollary \ref{main2}.

\bcorollary\label{va} If $(\mathscr E,h)$ is a uniformly RC-positive
vector bundle over a compact complex manifold $X$. Then \beq H^0(X,
(\mathscr E^*)^{\ts m})=0, \qtq{for every $m\geq 1$.}\eeq
\ecorollary \bproof It follows from Proposition \ref{tensor},
Corollary \ref{relation} and Theorem \ref{vanishing}. \eproof

\noindent The following special case is of particular interest.
\bcorollary Let $X$ be a compact complex manifold. If there exists a
smooth Hermitian metric $h$ such that $(T_X,h)$ is uniformly
RC-positive, then $H_{\bp}^{p,0}(X)=0$ for $1\leq p\leq \dim X$ and
\beq H^0(X,(T^*_X)^{\ts m})=0, \qtq{for all} m\geq 1.\eeq
\ecorollary \bproof By Corollary \ref{relation} and Proposition
\ref{tensor}, we know $\Lambda^p T_X$ and $T_X^{\ts m}$ are all
uniformly RC-positive. Hence by Corollary \ref{va},
$H_{\bp}^{p,0}(X)\cong H^0(X, \Lambda^pT_X^*)=0$. \eproof

\btheorem\label{asy} Let $(\mathscr E,h)$ be a uniformly RC-positive
vector bundle over a compact complex manifold $X$. Then for any line
bundle $\mathscr L$ over $X$, there exists a positive constant
$C=C(\mathscr L)$ such that $\mathscr E^{\ts m}\ts \mathscr L^{\ts
k}$ is uniformly RC-positive for all $m,k\in \N^+$ with $m\geq Ck$.
In particular,  \beq H^0\left(X,(\mathscr E^*)^{\ts m}\ts (\mathscr
L^*)^{\ts k}\right)=0.\label{rationalvanishing}\eeq

\etheorem

\bproof We fix an arbitrary smooth Hermitian metric $g$ on $\mathscr
L$, and assume that the curvature of $(\mathscr L,g)$ is bounded
from  by a negative constant $-B$. Fix a point $q\in X$. Since
$(\mathscr E,h)$ is uniformly RC-positive, there exists a positive
constant $C$ and a unit vector $u\in T_qX$ such that $$\inf_{v\in
\mathscr E_q\setminus\{0\}}\frac{R^\mathscr E(u,\bar u,v,\bar
v)}{|v|_h^2}\geq C>0$$
 By a similar computation as in Proposition \ref{tensor}, for any $v\in \mathscr E_q^{\ts
 m}$, we have
\beq R^{(\mathscr E^{\ts m},h^{\ts m})}(u,\bar u,v,\bar v)\geq
mC|v|^2.\eeq We choose a local unitary frame $e\in \Gamma(X,\mathscr
L)$ centered at $q$, then \be R^{\mathscr E^{\ts m}\ts \mathscr
L^{\ts k}}(u,\bar u,v\ts e^{\ts k},\bar{v\ts e^{\ts
k}})&=&R^{(\mathscr E^{\ts m},h^{\ts m})}(u,\bar u,v,\bar v)+k|v|^2
R^{\mathscr L}(u,\bar u)\\&\geq& mC|v|^2-kB|v|^2.\ee Hence if $m\geq
kB/C+1$, then
$$ R^{\mathscr E^{\ts m}\ts \mathscr L^{\ts k}}(u,\bar u,v\ts e^{\ts
k},\bar{v\ts e^{\ts k}})\geq |v|^2.$$ Hence, $\mathscr E^{\ts m}\ts
\mathscr L^{\ts k}$ is uniformly RC-positive. By Corollary \ref{va},
we obtain (\ref{rationalvanishing}). \eproof

\noindent Recall that $\{z^i\}, \{e^\alpha\}$ are the local
holomorphic coordinates and holomorphic frames on $X$ and $\mathscr
E$ respectively.

\blemma Let $\tilde h =e^{-f} h$ for some $f\in C^2(X,\R)$. Then the
curvature tensor $\tilde R$ of  $(\mathscr E,\tilde h)$ has the
expression \beq \tilde R_{i\bar j \alpha\bar \beta}=e^{-f}(R_{i\bar
j \alpha\bar\beta}+f_{i\bar j} h_{\alpha\bar
\beta}),\label{conformal}\eeq where $R$ is the curvature tensor of
$(\mathscr E,h)$. \elemma \bproof It follows by a
standard computation. 
 \eproof

\btheorem\label{perturb} Let $(\mathscr E,h)$ be a Hermitian
holomorphic vector bundle over
 a compact complex manifold $X$. Suppose $(\mathscr E,h)$ is uniformly RC-semi-positive over $X$. If there exist two  open
 subsets $S$ and $U$ of $X$ such that
\bd \item $U$ is strongly pseudoconvex and $\bar S\subset U$;

\item $(\mathscr E,h)$ is uniformly RC-positive on $X\setminus S$.

\ed Then there
 exists a smooth Hermitian metric $\tilde h$ on $\mathscr E$ such that $(\mathscr E,\tilde
 h)$ is uniformly RC-positive over $X$.
\etheorem \bproof Fix an arbitrary smooth Hermitian metric $\omega$
on $X$.  We define \beq C=\inf_{q\in X\setminus S} \sup_{u\in T_q
X\setminus\{0\}} \inf_{v\in \mathscr E_q\setminus\{0\}}
\frac{R(u,\bar u,v,\bar v)}{|u|^2_\omega|v|_h^2}.\label{C1}\eeq
Since $X\setminus S$ is compact and $(\mathscr E,h)$ is uniformly
RC-positive over $X\setminus S$, it is easy to see that $C>0$. There
exists a ``cut-off" function $f\in C^\infty(X,\R)$ such that

\bd
\item over $X$, we have
\beq  (\sq \p\bp f)(u,\bar u)\geq -\frac{C}{2}|u|_\omega^2;
\label{cut}\eeq
\item over $\bar S$, we have \beq (\sq
\p\bp f)(u,\bar u)\geq C_1 |u|_\omega^2\label{c3}\eeq for some
positive constant $C_1$.

\ed

\noindent Indeed, since $U$ is strongly pseudoconvex, there exists a
smooth strictly plurisubharmonic function $\phi \in
\mathrm{Psh}(U)$. In particular, there exists a positive constant
$\tilde C_1$ such that $ (\sq \p\bp \phi)(u,\bar u)\geq \tilde C_1
|u|_\omega^2$ over the compact set $\bar S$. Next, we can extend the
smooth function $\phi|_{\bar S}$ to $X$ and get a new function
$\Phi\in C^\infty(X).$ It is obvious that, there exists a positive
constant $\tilde C$ such that $(\sq \p\bp \Phi)(u,\bar u)\geq
-\tilde C|u|_\omega^2$ over $X$. Now we define $f=\frac{C}{2\tilde
C}\Phi$, then $f$ satisfies (\ref{cut}) and (\ref{c3}) with
$C_1=\frac{C \tilde C_1}{2\tilde C}$.

We define a new smooth Hermitian metric $\tilde h=e^{-f}h$ on
$\mathscr E$. By formula (\ref{conformal}), the curvature tensor
$\tilde R$ of $(\mathscr E,\tilde h)$ satisfies \beq \tilde R(u,\bar
u, v,\bar v)=e^{-f}\left(R(u,\bar u, v,\bar v)+(\p\bp f)(u,\bar
u)\cdot |v|_h^2\right).\eeq  We claim that $(\mathscr E,\tilde h)$
is uniformly RC-positive over $X$. Indeed, for a point $q\in S$,
since $(\mathscr E,h)$ is uniformly RC-semi-positive at $q$,  there
exists a unit vector $u\in T_q X$ such that $R(u,\bar u, v,\bar
v)\geq 0$. By estimates (\ref{c3}), we have  \be \tilde R(u,\bar u,
v,\bar v)&=&e^{-f}\left(R(u,\bar u, v,\bar v)+(\p\bp f)(u,\bar
u)\cdot |v|_h^2\right)\\&\geq& e^{-f}(\p\bp f)(u,\bar u)\cdot
|v|_h^2\\&\geq& C_1 e^{-f}|v|^2_h. \ee for all $v\in \mathscr E_q$.
For a point $q\in X\setminus S$, since $(\mathscr E,h)$ is uniformly
RC-positive over $X\setminus S$, by formula (\ref{C1}) and
Proposition \ref{est2}, there exists some unit vector $u\in T_q X$
such that \beq R(u,\bar u, v,\bar v)\geq C|v|_h^2\eeq for every
$v\in \mathscr E_q$. On the other hand, by estimate (\ref{cut}), we
have \beq \tilde R(u,\bar u, v,\bar v)=e^{-f}\left(R(u,\bar u,
v,\bar v)+(\p\bp f)(u,\bar u)\cdot |v|_h^2\right)\geq \frac{C}{2}
e^{-f} |v|^2_h. \eeq Hence, we conclude $(\mathscr E,\tilde h)$ is
uniformly RC-positive over $X$. \eproof

\noindent \emph{The proof of Theorem \ref{quasi2}.}
 By Theorem \ref{perturb}, $T_X$ has a   uniformly
RC-positive metric $\tilde h=e^{-f}\cdot h$.  By Proposition
\ref{tensor}, we know $\Lambda^2T_X$ is uniformly RC-positive.  By
Corollary \ref{relation} and Corollary \ref{kp}, $X$ is a projective
manifold. On the  other hand, by Theorem \ref{asy}, for any line
bundle $\mathscr L$ over $X$, there exists a positive constant
$C=C(\mathscr L)$ such that  \beq H^0\left(X,(T_X^*)^{\ts m}\ts
\mathscr L^{\ts k}\right)=0\label{C5}\eeq for all $m,k\in \N^+$ with
$m\geq Ck$. Therefore, by a celebrated theorem of
Campana-Demailly-Peternell \cite[Theorem~1.1]{CDP},  $X$ is
rationally connected. Indeed, we deduce from (\ref{C5}) that for
each $1\leq p\leq \dim X$,  any invertible sheaf $\mathcal{F}\subset
\Om^p_X$ is not pseudo-effective. Otherwise, if there exists a
pseudo-effective invertible sheaf $\mathcal{F}\subset \Om^p_X$, then
there exists a very ample line bundle $A$ such that
$H^0(X,\mathcal{F}^{\ts \ell}\ts A)\neq 0$ for all $\ell\geq 0$, and
so
 $H^0(X,\mathrm{Sym}^{\ts \ell}\Om_X^p\ts A)\neq 0$. Since for some large $m$, $\mathrm{Sym}^{\ts \ell}\Om_X^p\subset (T_X^{*})^{\ts m}$, we get a  contradiction
 to (\ref{C5}).  In particular, when $p=n$, $K_X$ is not
pseudo-effective. Thanks to \cite{BDPP13}, $X$ is uniruled. Let
$\pi:X\dashrightarrow Z$ be the associated MRC fibration of $X$.
After possibly resolving the singularities of $\pi$ and $Z$, we may
assume that $\pi$ is a proper morphism and $Z$ is smooth. By a
result of Graber, Harris and Starr \cite[Corollary~1.4]{GHS03}, it
follows that the target $Z$ of the MRC fibration is either a point
or a positive dimensional variety which is not unruled. Suppose $X$
is not rationally connected, then $\dim Z\geq 1$. Hence $Z$ is not
uniruled, by \cite{BDPP13} again, $K_Z$ is pseudo-effective. Since
 $K_Z=\Om^{\dim Z}_Z\subset
\Om_X^{\dim Z}$ is pseudo-effective, we get a contradiction. Hence
$X$ is rationally connected. \qed

\vskip 1\baselineskip

\noindent If  $S=\varnothing$, we have \bcorollary\label{urcrc} Let
$X$ be a compact K\"ahler manifold. Suppose $X$ has a Hermitian
metric $h$ such that $(T_X,h)$ is uniformly RC-positive, then $X$ is
a projective and rationally connected manifold. \ecorollary

\bremark Corollary \ref{urcrc}  can also be proved by using
Proposition \ref{tensor} and \cite[Theorem~1.3]{Yang18}. The proof
here is slightly simpler. \eremark

\noindent By using rational connectedness, one has \bcorollary  Let
$X$ be a compact K\"ahler manifold. If $X$ admits a smooth Hermitian
metric $\omega$ such that $(T_X,\omega)$ is uniformly RC-positive,
then there is  no non-constant holomorphic map from $X$ to a compact
complex manifold $Y$ where $Y$ lies in one of the following \bd\item
$Y$ is Kobayashi hyperbolic;
\item $Y$ has nef cotangent bundle; \item  $Y$ has a Hermitian metric
with non-positive holomorphic sectional curvature;

\item $Y$ contains no rational curve.

\ed \ecorollary

\vskip 2\baselineskip

\section{RC-positive metrics on rationally connected manifolds}
\label{re}

 In this section, we will discuss general theory for uniformly RC-positive vector bundles and prove
 Theorem \ref{main3}. Let's recall that a line bundle $\mathscr  L$ is uniformly RC-positive if and
 only if it has a smooth Hermitian metric $h$ such that its
 curvature $-\sq\p\bp\log h$ has at least one positive eigenvalue
 everywhere. In \cite[Theorem~1.4]{Yang17D}, we obtained an
 equivalent characterization for uniformly RC-positive line
 bundles.
 \btheorem\label{equi} Let $\mathscr L$ be a holomorphic line bundle over a compact complex manifold $X$. The
 following statements  are equivalent.

 \bd \item $\mathscr L$ is uniformly RC-positive;
 \item the dual line bundle $\mathscr L^{*}$ is not pseudo-effective.
 \ed
 \etheorem
\noindent Two key ingredients in the proof of Theorem \ref{equi} are
a conformal (exponential) perturbation method  and an  integration
criterion for pseudo-effectiveness over compact complex manifolds(
e.g. \cite{La99}). We refer to \cite{Yang17a, Yang17D} and the
references therein.

\bcorollary Let $X$ be a projective manifold. If $\mathscr L$ is
uniformly RC-semi-positive, then the dual line bundle $\mathscr
L^{*}$ is not big. \bproof Let $\mathscr A$ be an ample line bundle
over $X$. We argue by contradiction. Suppose $\mathscr L^*$ is big,
then there exists a large number $k\in \Z_+$ such that $\mathscr
L^{*k}\ts \mathscr A^{*}$ is big. By Theorem \ref{equi}, $\mathscr
L^k\ts \mathscr A$ can not be uniformly RC-positive which is absurd.
\eproof \ecorollary \noindent

\noindent The following concepts are generalizations of uniformly
RC-positivity for line bundles.
\begin{definition} Let $\mathscr L$ be a  line bundle
over a  compact complex manifold $X$. \bd \item $\mathscr L$ is
called \emph{$q$-positive}, if there exists a smooth Hermitian
metric $h$ on $\mathscr L$ such that the Chern curvature
$R^{(\mathscr L,h)}=-\sq\p\bp\log h$ has at least $(\dim X-q)$
positive eigenvalues at every point on $X$.
\item $\mathscr L$ is called \textit{ $q$-ample}, if for any coherent sheaf
$\mathscr{F}$ on $X$ there exists a positive integer
$m_{0}=m_{0}(X,\mathscr L,\mathscr{F})>0$ such that
$$ H^{i}\left( X, \mathscr{F}\otimes \mathscr L^{ m}\right) =0\ \ \
\mathrm{for}\  i>q,\ m\geq m_{0}. $$\ed
\end{definition}
\noindent Note that, \emph{the nontrivial range is $0\leq q\leq \dim
X-1$.} By the celebrated Cartan--Serre--Grothendieck criterion for
ampleness and Kodaira embedding theorem, one has
\btheorem[Cartan--Serre--Grothendieck, Kodaira]\label{CSGK} Let
$\mathscr L\>X$ be a holomorphic line bundle over a projective
manifold $X$. Then the following statements are equivalent \bd \item
$\mathscr L$ is $0$-ample;
\item $\mathscr L$ is $0$-positive. \ed\etheorem

\noindent  As a weak dual to the Cartan--Serre--Grothendieck-Kodaira
Theorem \ref{CSGK}, we established in \cite{Yang17D} that

\btheorem\label{AG} Let $\mathscr L\>X$ be a holomorphic line bundle
over a projective manifold $X$. Then the following statements are
equivalent \bd \item $\mathscr L$ is $(\dim X-1)$-ample;
\item $\mathscr L$ is $(\dim X-1)$-positive;
  \item $\mathscr L$ is uniformly
 RC-positive. \ed\etheorem

\noindent For related results on $q$-ampleness and $q$-positivity,
we refer to \cite{AG62, DPS93, Dem11, Tot13, Ott12, Bro12, Mat13,
Yang17D} and the references therein.\\

  Let $\mathscr E$ be a vector bundle over $X$ and $\P(\mathscr E^*)$ be the projective
bundle of $\mathscr E$. The points of the projective bundle
$\P(\mathscr E^*)$ of $\mathscr E\>X$ can be identified with the
hyperplanes of $\mathscr E$. Note that every hyperplane $V$ in
$\mathscr E_z$ corresponds bijectively to the line of linear forms
in $\mathscr E^*_z$ which vanish on $V$.
 Let $\pi : \P(\mathscr E^*) \> X$ be
the natural projection. There is a tautological hyperplane subbundle
$S$ of $\pi^*\mathscr E$ over $\P(\mathscr E^*)$ such that
$$S_{[\xi]} = \xi^{-1}(0) \subset \mathscr E_z$$ for
all $\xi\in \mathscr E_z^*\setminus\{0\}$. The quotient line bundle
$\pi^*\mathscr E/S$ is denoted $\sO_\mathscr E(1)$ and is called the
tautological line bundle associated to $\mathscr E\>X$. Hence there
is an exact sequence of vector bundles over $\P(\mathscr E^*)$ \beq
0 \>S\>\pi^*\mathscr E\>\sO_\mathscr E(1)\> 0.\eeq A holomorphic
vector bundle $\mathscr E\>X$ is called ample if the line bundle
$\sO_\mathscr E(1)$ is ample over $\P(\mathscr E^*)$.

  \bproposition\label{pro} Let $(\mathscr E,h)$ be  RC-positive vector bundle over a compact complex manifold $X$. Then \bd \item the
  tautological line bundle $\sO_\mathscr E(1)$ is $(\dim X-1)$-positive over
  $\P(\mathscr E^*)$;

\item  $\sO_{\mathscr E^*}(-1)$ is uniformly RC-positive over
  $\P(\mathscr E)$;
\item  $\sO_{\mathscr E^*}(1)$ is not pseudo-effective over
  $\P(\mathscr E)$.
  \ed
\eproposition

\bproof It follows from the curvature formulas of $\sO_{\mathscr
E^*}(-1)$ and $\sO_\mathscr E(1)$ induced by $(\mathscr E, h)$ (e.g.
\cite[Proposition~4.1]{Yang18}, \cite[Proposition~2.6]{Yang18a}) and
Theorem \ref{equi}. \eproof


\btheorem\label{Theorem} Let $\mathscr E$ be a holomorphic vector
bundle over a projective manifold $X$. Then the following statements
are equivalent:

\bd \item $\sO_\mathscr E(1)$ is $(\dim X-1)$-ample;

\item $\sO_{\mathscr E^*}(1)$ is not pseudo-effective;

\item $\sO_{\mathscr E^*}(-1)$ is uniformly RC-positive.

\ed

\etheorem

\bproof If $\sO_\mathscr E(1)$ is $(\dim X-1)$-ample, we deduce that
$\sO_{\mathscr E^*}(1)$ is not pseudo-effective. Otherwise, if
$\sO_{\mathscr E^*}(1)$ is pseudo-effective, it is well-known that
there exists an ample line bundle $\mathscr L$ on $\P(\mathscr E)$
such that \beq H^0(\P(\mathscr E),\sO_{\mathscr E^*}(m)\ts \mathscr
L)\neq 0 \qtq{for all $m\geq 0$.} \label{non}\eeq  More precisely,
 $\mathscr L$ can be chosen in such a way: fix a very ample line bundle $\mathscr H$ over
$\P(\mathscr E)$, and if $\mathscr L$ is an ample line bundle such
that
$$\mathscr L\ts K_{\P(\mathscr E)}^{-1}\ts \mathscr H^{-(\dim \P(\mathscr E)+1)}$$ is ample, then
(\ref{non}) holds.  We can choose $\mathscr L$ in a special form
\beq \mathscr L=\sO_{\mathscr E^*}(m_0)\ts \pi^*(\mathscr A^{\ts
m_1})\label{pre}\eeq where $\mathscr A$ is an ample line bundle over
$X$,
 and $m_0,m_1$ are two large positive integers with $m_1\gg m_0$ and $\pi:\P(\mathscr E)\>X$ is the
 projection. For reader's convenience, we write down the construction explicitly. Indeed, for large $k$, $\mathscr E^*\ts \mathscr A^{\ts k}$ is an ample
 vector bundle over $X$, i.e. $\sO_{\mathscr E^*}(1)\ts \pi^*(\mathscr A^{\ts k})$ is
 an ample line bundle over $\P(\mathscr E)$. There exists a large $k_1$ such
 that $\mathscr H=\sO_{\mathscr E^*}(k_1)\ts \pi^*(\mathscr A^{\ts (kk_1)})$ is
 a very ample line bundle over $\P(\mathscr E)$. On the other hand,
 $$K_{\P(\mathscr E)}=\sO_{\mathscr E^*}(-r)\ts \pi^*(K_X)\ts \pi^*(\det \mathscr E^*)$$ where
 $r$ is the rank of $\mathscr E$.
There exists a positive number $k_2$ such that $$\mathscr A^{\ts
k_2}\ts K_X^{-1}\ts \det \mathscr E$$ is ample over $X$. Now we can
take
$$m_0=k_1(n+r)-r+1 \qtq{and} m_1=k_2+kk_1(n+r)+k$$
and so $$\mathscr L\ts K_{\P(\mathscr E)}^{-1}\ts \mathscr
H^{-(n+r)}=\sO_{\mathscr E^*}(1)\ts \pi^*(\mathscr A^{\ts k})\ts
\pi^*\left(\mathscr A^{\ts k_2}\ts K_X^{-1}\ts \det \mathscr
E\right)$$ is ample over $\P(\mathscr E)$.

 Therefore, by the Le Potier isomorphism,  (\ref{non}) is equivalent
 to
\beq H^0(X,\mathrm{Sym}^{\ts k}\mathscr E^*\ts \mathscr A^{\ts
m_1})\neq 0 \label{nonva}\eeq for large $k$. Applying the Serre
duality on $X$ and the  Le Potier isomorphism again, we obtain
$$H^n(\P(\mathscr E^*), \sO_\mathscr E(k)\ts \pi_1^*((\mathscr A^*)^{\ts m_1})\ts \Om^n_{\P(\mathscr E^*)})\neq 0,$$
for large $k$ where $\pi_1:\P(\mathscr E^*)\>X$ is the projection.
Let $\mathscr{F}=\pi_1^*((\mathscr A^*)^{\ts m_1})\ts
\Om^n_{\P(\mathscr E^*)}$, and we know $\sO_\mathscr E(1)$ can not
be $(n-1)$-ample. This is a contradiction. Hence $\sO_{\mathscr
E^*}(1)$ is not pseudo-effective.  The proof of $(2)\Longrightarrow
(1)$ is similar. The equivalence of $(2)$ and $(3)$ follows from
Theorem \ref{equi}.
 \eproof

\noindent On the other hand, we have \bproposition \label{F} Let $X$
be a projective manifold. Suppose $\sO_\mathscr E(1)$ is $(\dim
X-1)$-ample, then
 the restriction of $\sO_\mathscr E(1)$ to every  smooth submanifold $Y\subset \P(\mathscr E^*)$ with  $\dim Y=\dim X$ is $(\dim X-1)$-positive.
\eproposition \bproof Let $f:Y\>\P(\mathscr E^*)$ be the inclusion
map. Using the projection formula and the Leray spectral sequence,
one has
$$H^i(Y, \mathscr F\ts (f^*\sO_\mathscr E(m)))=H^i(\P(\mathscr E^*), f_*(\mathscr F)\ts \sO_\mathscr E(m))$$
Hence, if $\sO_\mathscr E(1)\>\P(\mathscr E^*)$ is $(\dim
X-1)$-ample, $f^*(\sO_\mathscr E(1))\>Y$ is also $(\dim X-1)$-ample.
On the other hand, since $\dim Y=\dim X$, by Theorem \ref{AG}, the
$(\dim X-1)$-ample line bundle $f^*(\sO_\mathscr E(1))$ over $Y$ is
$(\dim X-1)$-positive. \eproof

\noindent As motivated by these properties, we propose the following
conjecture.

\begin{conjecture}\label{conjecture} Let $\mathscr E$ be a holomorphic vector bundle over a projective manifold $X$. Then the following statements are equivalent. \bd
\item $\sO_{\mathscr E^*}(-1)$ is RC-positive;
   \item $\sO_\mathscr E(1)$ is  $(\dim X-1)$-ample;
\item $\sO_\mathscr E(1)$ is $(\dim X-1)$-positive;
\item $\mathscr E$ is  RC-positive.
\ed
\end{conjecture}

\noindent Note that, the implications $(4)\Longrightarrow
(3)\Longrightarrow (2)\Longleftrightarrow (1)$ in Conjecture
\ref{conjecture} are known by Proposition \ref{pro}, \cite{AG62} and
 Theorem
\ref{Theorem}. On the other hand, when $\mathrm{rank}(\mathscr E)=1$
or $\dim X=1$, Conjecture \ref{conjecture} is
 true by Theorem \ref{AG} and \cite{CF90}. Note also that Conjecture
 \ref{conjecture} is also analogous to  a conjecture of Griffiths.

 \begin{conjecture} Let $\mathscr E$ be a vector bundle over a projective manifold $X$. Then the following statements are equivalent.
\bd \item $\sO_\mathscr E(1)$ is positive;
\item $\mathscr E$ is Griffiths positive.
\ed
 \end{conjecture}

\noindent It is easy to see that if
$\mathscr E$ is Griffiths positive, then so is $\sO_\mathscr E(1)$.\\

As an application of the vector bundle theory discussed above, we
obtain a differential geometric characterization of rationally
connected manifolds, and  Theorem \ref{main3} is also a
 special case of it.

\btheorem\label{BBBQ} Let $X$ be a projective manifold. Then the
following statements are equivalent \bd \item $X$ is rationally
connected;
\item  the line bundle $\sO_{\Lambda^pT_X^*}(-1)$ is uniformly
RC-positive for every $1\leq p\leq \dim X$. \ed\etheorem

\bproof If $X$ is rationally connected, then by a variant of
  \cite[Theorem~1.1]{CDP} (e.g. \cite[Proposition~1.4]{Cam16}), for any ample line
bundle $\mathscr L$ on $X$, there exists a positive integer $m_0$
such that
$$ H^0(X,\mathrm{Sym}^{\ts m}\left(\Lambda^pT_X^*\right)\ts \mathscr L^{\ts k})=0
$$ for $m\geq m_0k$ and all $1\leq p\leq \dim X$. We claim $\sO_{\Lambda^pT_X^*}(-1)$ is RC-positive.
Otherwise, by Theorem \ref{Theorem}, $\sO_{\Lambda^pT_X^*}(1)$ is
pseudo-effective. Hence,  by using a similar proof as  in Theorem
\ref{Theorem},  we can find an ample line bundle $\mathscr A$ over
$X$ such that (\ref{nonva}) holds for $\mathscr E=\Lambda^p T_X$,
that is
$$H^0(X,\mathrm{Sym}^{\ts m}\left(\Lambda^pT_X^*\right)\ts \mathscr A^{\ts m_1})\neq 0$$
for all large $m$. This is a contradiction.

 On the other hand, if
$\sO_{\Lambda^pT_X^*}(-1)$ is uniformly RC-positive over
$\P(\Lambda^p T_X)$, then by the Le Potier isomorphism and  Theorem
\ref{asy}, for any coherent sheaf of the  form
$\mathscr{F}=\pi^*(\mathscr L^{\ts k})$ over $\P(\Lambda^p T_X)$, we
have \beq H^0(X,\mathrm{Sym}^{\ts m}\left(\Lambda^pT_X^*\right)\ts
\mathscr L^{\ts k} )\cong H^0(\P(\Lambda^p
T_X),\sO_{\Lambda^pT_X^*}(m)\ts \mathscr{F} )=0 \eeq for large $m$.
 Therefore,
by \cite[Theorem~1.1]{CDP} or \cite[Proposition~1.4]{Cam16}, $X$ is
rationally connected. \eproof

\noindent By Theorem \ref{BBBQ}, the following conjecture is a
special case of Conjecture \ref{conjecture}.
\begin{conjecture}\label{Q1} Let $X$ be a projective manifold.  If  $X$ is  rationally connected, then $T_X$ is  RC-positive.
\end{conjecture}

\noindent More generally,

\begin{problem} The following statements are equivalent on a projective manifold $X$.
\bd\item   $X$ is rationally connected;

\item  $T_X$ is RC-positive;

\item $T_X$ is uniformly RC-positive.
\ed
\end{problem}

\noindent In \cite[Corollary~1.5]{Yang18}, we obtain the following
result. \bcorollary\label{BBQ} Let $X$ be a compact K\"ahler
manifold. If there exist a Hermitian metric $\omega$ on $X$ and a
(possibly different) Hermitian metric $h$ on $T_X$ such that \beq
\mathrm{tr}_{\omega}R^{(T_X,h)} \in \Gamma(X, \mathrm{End}(T_X))\eeq
is positive definite, then $X$ is projective and rationally
connected. \ecorollary

\noindent In particular, by
 the celebrated Calabi-Yau theorem (\cite{Yau78}), one gets the
classical result of Campana(\cite{Cam92}) and Koll\'ar-Miyaoka-Mori
(\cite{KMM92}) that Fano manifolds are rationally connected. We
propose a conjecture converse to Corollary \ref{BBQ}, which is also
analogous to the classical fact that a compact complex manifold is
Fano if and only if it has a Hermitian metric with positive
Chern-Ricci curvature.

\begin{problem}\label{Q3} The following statements are equivalent on a projective manifold $X$.  \bd \item $X$ is rationally connected;
\item there exist a \emph{Hermitian metric}
$\omega$ on $X$ and a (possibly different) Hermitian metric $h$ on
$T_X$ such that $ \mathrm{tr}_\omega R^{(T_X,h)}$ is positive
definite. \ed
\end{problem}
\bremark Problem \ref{Q3} is also known to J.-P. Demailly \cite{DC}.
A positive solution to Problem \ref{Q3} gives an affirmative  answer
to Conjecture \ref{Q1}. \eremark

\vskip 2\baselineskip

\section{Compact K\"ahler manifolds with nonnegative holomorphic sectional curvature}\label{HSC00}

A compact K\"ahler manifold  $(X,\omega)$ has positive (resp.
nonnegative) holomorphic sectional curvature, if for any nonzero
vector $\xi=(\xi^1,\cdots, \xi^n)$,
$$R_{i\bar j k\bar \ell}\xi^i\bar\xi^j\xi^k\bar\xi^\ell>0\ \ \text{(resp. $\geq 0$)} $$
at each point of $X$. The negativity and non-positivity of the
holomorphic sectional curvature can be defined in a similar way.

\btheorem\label{phscurc} Let $(X,\omega)$ be a compact K\"ahler
manifold with positive (resp.  nonnegative) holomorphic sectional
curvature, then $(T_X,\omega)$ is uniformly RC-positive (resp.
uniformly RC-semi-positive).

 \etheorem\bproof Let
$\kappa$ be the positive holomorphic sectional curvature of
$(X,\omega)$, i.e. \beq \kappa=\inf_{q\in X}\inf_{U\in T_q
X\setminus\{0\}}\frac{R(U,\bar U, U,\bar U )}{|U|^4}.\eeq For any
point $q\in X$, let $e\in T_q X$ be a unit vector such that \beq
R(e,\bar e, e,\bar e)=\inf_{U\in T_q X\setminus\{0\}}\frac{R(U,\bar
U, U,\bar U )}{|U|^4}.\eeq We know $R(e,\bar e, e,\bar e)\geq
\kappa> 0.$ On the other hand, by \cite[Lemma~6.1]{Yang18}, for any
unit vector $W\in T_q X$, we have \beq 2R(e,\bar e, W,\bar W)\geq
(1+|\la W, e\ra|^2)R(e,\bar e, e, \bar e)\geq \kappa.\eeq Hence, for
any vector $v\in T_qX$, we obtain
$$R(e,\bar e,
v,\bar v)\geq \frac{\kappa}{2}|v|^2,$$ and so $(T_X,\omega)$ is
uniformly RC-positive.
 \eproof

\bremark\label{remark} If a compact K\"ahler manifold $(X,\omega)$
hash negative (resp. nonpositive) holomorphic sectional curvature,
then $(T_X,\omega)$ is uniformly RC-negative (resp. uniformly
RC-nonpositive).  There are many  K\"ahler and non-K\"ahler complex
manifolds which have uniformly RC-positive tangent bundles, for
instances, \bd

 \item[$\bullet$] compact K\"ahler manifold with non-negative
holomorphic bisectional curvature and positive first Chern
class(\cite{Mok88});
\item[$\bullet$]  Hopf manifold $\S^1\times \S^{2n+1}$ (\cite[formula
(6.4)]{LY14}).

\ed

\eremark

\noindent As an application of Theorem \ref{phscurc} and Theorem
\ref{quasi2}, we obtain

  \btheorem \label{QHSC1} Let $(X,\omega)$ be a compact
K\"ahler manifold with  nonnegative holomorphic sectional curvature.
If there exist two open
 subsets $S$ and $U$ of $X$ such that $U$ is strongly pseudoconvex, $\bar S\subset U$ and $(X,\omega)$ has positive holomorphic sectional curvature on $X\setminus S$.
Then $X$ has a uniformly RC-positive Hermitian metric. In
particular, $X$ is a projective and rationally connected manifold.
 \etheorem

\bremark It is not hard  to see that Theorem \ref{QHSC1} can also
hold under certain weaker conditions.  For related topics on
 holomorphic sectional curvature, we refer to \cite{HW15,
ACH15, Yang16, Liu16, AHZ16, YZ16, AH17, Mat18} and the references
therein. \eremark

\vskip 2\baselineskip

\section{RC-positive Finsler vector bundles over  complex manifolds }\label{3}

Let $\mathscr E$ be a holomorphic vector bundle over a  complex
manifold $X$ with complex rank $r$.  Let $z^1,\cdots, z^n$ be the
local holomorphic coordinates on $X$ and $w^1,\cdots ,w^r$ be the
holomorphic coordinates on the fiber of $\mathscr E$. Let $\mathscr
E\setminus\{0\}$ be the complement of the zero section of $\mathscr
E$.

\bdefinition A pseudoconvex complex Finsler metric $\mathfrak{F}$ on
$\mathscr E$ is a continuous function $\mathfrak{F}: \mathscr
E\>[0,+\infty)$ satisfying \bd
\item $\mathfrak{F}$ is smooth on $\mathscr E\setminus\{0\}$; \item
$\mathfrak{F}(z,w)> 0$ for all $w\neq 0$; \item
$\mathfrak{F}(z,\lambda w)=|\lambda|^2 \mathfrak{F}(z,w)$ for all
$\lambda\in \C$;

\item The $(r\times r)$ Hermitian matrix  $\left(\frac{\p^2 \mathfrak F}{\p w^\alpha\p\bar
w^\beta}\right)$ is positive definite over $\mathscr E\setminus
\{0\}$.

 \ed \edefinition

\noindent Let $\mathfrak{F}$ be a pseudoconvex complex Finsler
metric on $\mathscr E\>X$. It is well-known that
$$\left(h_{\alpha\bar\beta}\right)=\left(\frac{\p^2 \mathfrak F}{\p w^\alpha\p\bar
w^\beta}\right)$$ defines a smooth Hermitian metric on the
holomorphic vector bundle $\pi^*\mathscr E\>\P(\mathscr E^*)$ where
$\pi:\P(\mathscr E^*)\>X$ is the projection. Note that, in general
$\mathfrak{F}$ does not give a Hermitian metric on $\mathscr E\>X$.

\bdefinition Let $\mathfrak{F}$ be a pseudoconvex complex Finsler
metric on $\mathscr E\>X$. $(\mathscr E,\mathfrak{F})$ is called an
RC-positive (resp. a uniformly RC-positive) Finsler vector bundle if
the induced Hermitian vector bundle $(\pi^*\mathscr E, h)$ is
RC-positive (resp. uniformly RC-positive) over $\P(\mathscr E^*)$.

\edefinition

\bremark If $(\mathscr E,h)$ is a Hermitian holomorphic vector
bundle, then it induces a pseudoconvex Finsler metric $\mathfrak{F}$
on $\mathscr E$. Moreover, if $(\mathscr E,h)$ is RC-positive (resp.
uniformly RC-positive), then $(\mathscr E,\mathfrak{F})$ is
RC-positive (resp. uniformly RC-positive). \eremark

\btheorem\label{F1} Let $X$ be a compact complex manifold. Suppose
$(\mathscr E,\mathfrak{F})$ is an RC-positive Finsler vector bundle.
Then we have
$$H^0(X, \mathrm{Sym}^{\ts m}\mathscr E^*)=0, \qtq{for all} m\geq 1.$$
Moreover, for any line bundle $\mathscr \mathscr \mathscr L\>X$,
there exists a positive constant $c_\mathscr L$ such that
$$H^0(X, \mathrm{Sym}^{\ts m}\mathscr E^*\ts \mathscr  L^{*\ts k})=0,$$
for all positive integers $m,k$ with $m\geq c_\mathscr Lk$.
 \etheorem

\bproof If $(\mathscr E,\mathfrak{F})$ is an RC-positive Finsler
bundle, then the induced Hermitian vector bundle $(\pi^*\mathscr E,
h)$ is RC-positive. Let $\pi:\P(\mathscr E^*)\>X$ be the natural
projection. Let $\mathscr F=\pi^*(\mathscr E)$ and $Y=\P(\mathscr
E^*)$. By Proposition \ref{pro}, $\sO_{\mathscr F^*}(-1)$ is an
RC-positive line bundle over the projective bundle $\tilde \pi:
\P(\mathscr F)\>Y$. Hence, by Theorem \ref{asy}, we have \beq
H^0\left(\P(\mathscr F), \sO_{\mathscr F^*}(m)\ts
\tilde\pi^*\left(\pi^*(\mathscr L^{*\ts k})\right)\right)=0, \eeq
for all positive integers $m,k$ with $m\geq c_\mathscr Lk$. By the
Le Potier isomorphism, we have  \be H^0(X, \mathrm{Sym}^{\ts
m}\mathscr E^*\ts \mathscr L^{\ts k})&\cong& H^0(\P(\mathscr
E^*),\pi^*\left(\mathrm{Sym}^{\ts m}\mathscr E^*\right)\ts
\pi^*(\mathscr L^{\ts k}))\\&\cong& H^0(\P(\mathscr
E^*),\mathrm{Sym}^{\ts m}\pi^*\left(\mathscr E^*\right)\ts
\pi^*(\mathscr L^{\ts k}))\\&=&
 H^0(Y, \mathrm{Sym}^{\ts
m}\mathscr F^*\ts \pi^*(\mathscr L^{*\ts k}))\\&\cong&
H^0\left(\P(\mathscr F), \sO_{\mathscr F^*}(m)\ts
\tilde\pi^*\left(\pi^*(\mathscr L^{*\ts k})\right)\right)\\&=&0. \ee
If we take $\mathscr L$ to be a trivial line bundle, then there
exists a large positive integer $m$ such that $H^0(X,
\mathrm{Sym}^{\ts m}\mathscr E^*)=0.$ It is easy to see that $H^0(X,
\mathscr E^*)=0$ and so $H^0(X, \mathrm{Sym}^{\ts m}\mathscr E^*)=0$
for all $m\geq 1$. \eproof

\btheorem\label{F2} Let $X$ be a compact complex manifold. Suppose
$(\mathscr E,\mathfrak{F})$ is a uniformly RC-positive Finsler
vector bundle. Then we have
$$H^0(X, (\mathscr E^*)^{\ts m})=0, \qtq{for all} m\geq 1.$$
Moreover, for any line bundle $\mathscr L\>X$, there exists a
positive constant $c_\mathscr L$ such that
$$H^0(X, (\mathscr E^*)^{\ts m}\ts \mathscr L^{\ts k})=0,$$
for all positive integers $m,k$ with $m\geq c_\mathscr Lk$.
 \etheorem

\bproof The induced Hermitian vector bundle $(\pi^*\mathscr E, h)$
is uniformly RC-positive. By Proposition \ref{tensor},
$\pi^*(\mathscr E^{\ts m})\cong(\pi^*(\mathscr E))^{\ts m}$ is
uniformly RC-positive for all $m\geq 1$. On the other hand,
$$H^0\left(X, (\mathscr E^*)^{\ts m}\ts \mathscr L^{\ts k}\right)\cong H^0\left(\P(\mathscr E^*), (\pi^*(\mathscr E^*))^{\ts m}\ts \pi^*(\mathscr L^{\ts k})\right).$$
Hence, Theorem \ref{F2} follows from Theorem \ref{asy}. \eproof

\noindent As an application of Theorem \ref{F1} and
\cite[Theorem~1.4]{Yang18}, we obtain

\btheorem\label{F3} Let $X$ be a compact K\"ahler manifold of
complex dimension $n$. Suppose that for every $1\leq p\leq n$, there
exists a pseudoconvex Finsler metric $\mathfrak{F}_p$ on  $\Lambda^p
T_X$ such that $(\Lambda^p T_X, \mathfrak{F}_p)$ is RC-positive,
then $X$ is projective and rationally connected. \etheorem

\noindent Similarly, as an application of Theorem \ref{F2} and
Corollary \ref{main2}, we obtain

\btheorem \label{main20} Let $X$ be a compact K\"ahler manifold. If
$X$ admits a pseduoconvex Finsler metric $\mathfrak{F}$ such that
$(T_X,\mathfrak{F})$ is uniformly RC-positive, then $X$ is
projective and rationally connected. \etheorem

\bremark It is well-known that there exists a one-to-one
correspondence between the set of Hermitian metrics on
$\sO_{\mathscr E^*}(-1)$ and the set of Finsler metrics on $\mathscr
E^*$ or $\mathscr E$. We can also define a Finsler vector bundle
$(\mathscr E, \mathfrak{F})$ to be RC-positive, if the induced
Hermitian metric on $\sO_{\mathscr E^*}(-1)$ is RC-positive.
\eremark

\vskip 2\baselineskip


\begin{thebibliography}{kawamata}




\bibitem[ACH15]{ACH15}
Alvarez, A.; Chaturvedi, A.; Heier, G.
 Optimal pinching for the holomorphic sectional curvature of {H}itchin's
  metrics on Hirzebruch surfaces.
 \emph{Contemp. Math.}, 133--142, 2015.

\bibitem[AG62]{AG62}
Andreotti, A.; Grauert, H. \textit {Th\'eor\`eme de finitude pour la
cohomologie des espaces complexes.} \emph{Bull. Soc. Math. France}
{\bf 90} (1962), 193--259. 

\bibitem[AH17]{AH17} Chaturvedi, A.; Heier, G. Hermitian metrics of positive holomorphic sectional curvature on
fibrations.
\href{https://arxiv.org/abs/1707.03425}{arXiv:1707.03425v1}

\bibitem[AHZ16]{AHZ16}
Alvarez, A.; Heier, G.; Zheng, F.-Y.
 On projectivized vector bundles and positive holomorphic sectional
  curvature.  Proc. Amer. Math. Soc. \textbf{146} (2018),  2877--2882. 



\bibitem[Bro12]{Bro12} Brown, M.-V. Big $q$-ample line bundles. Compos. Math.
\textbf{148} (2012), no. 3, 790--798. 


\bibitem[BC15]{BC15} Brunebarbe, Y.; Campana, F. Fundamental group and pluri-differentials on compact K\"ahler
manifolds.  \emph{Mosc. Math. J.}  \textbf{16} (2016), no. 4,
651--658. 

\bibitem[BDPP13]{BDPP13} Boucksom, S.;  Demailly, J.-P.;  Paun, M.;  Peternell, P. The
pseudoeffective cone of a compact K\"ahler manifold and varieties of
negative Kodaira dimension. \emph{J. Algebraic Geom.} \textbf{22}
(2013) 201--248. 







\bibitem[Cam92]{Cam92} Campana, F. {Connexit\'e rationnelle des vari\'et\'es de Fano.}  \emph{Ann. Sci. Ecole Norm. Sup.} (4) \textbf{25} (1992), no. 5, 539--545.
\bibitem[Cam16]{Cam16} Campana, F. {Slope rational connectedness for orbifolds}. \href{https://arxiv.org/abs/1607.07829}{arXiv:1607.07829}.






\bibitem[CDP14]{CDP}  Campana, F.; Demailly, J.-P.; Peternell, T. {Rationally connected manifolds and semipositivity of the Ricci curvature.} in Recent advances in Algebraic Geometry. \emph{LMS Lecture Notes Series}  417. (2014),
71--91. 

\bibitem[CF90]{CF90}   Campana, F.;  Flenner, H. {A characterization of ample vector bundles on a curve}.\emph{ Math. Ann.} \textbf{287} (1990), no. 4,
571--575. 

\bibitem[CP14]{CP14} Campana, F.;  Paun, M.
{Positivity properties of the bundle of logarithmic tensors on
compact K\"ahler manifolds.} \emph{Compositio Math.} \textbf{152}
(2016) 2350--2370.

\bibitem[CH17]{CH17} Cao, J.-Y.; Hoering, A.. A decomposition theorem for projective manifolds with nef anticanonical
bundle. \href{https://arxiv.org/abs/1706.08814}{arXiv:1706.08814}

\bibitem[Dem11]{Dem11}
Demailly, J.-P. {A converse to the Andreotti-Grauert theorem.} Ann.
Fac. Sci. Toulouse Math. (6) {\bf 20} (2011), Fascicule Special,
123--135. 

\bibitem[Dem]{DC} Demailly, J.-P.  Private communication.

\bibitem[DPS93]{DPS93}
Demailly, J.-P.; Peternell, T.; Schneider, M. {Holomorphic line
bundles with partially vanishing cohomology.} \emph{Proceedings of
the Hirzebruch 65 Conference on Algebraic Geometry} (Ramat Gan,
1993), 165--198, Israel Math. Conf. Proc., {\bf 9}, Bar-Ilan Univ
1996. 

\bibitem[DPS96B]{DPS96B} Demailly, J.-P.; Peternell, T.; Schneider, M. Compact K\"ahler manifolds with hermitian semipositive anticanonical
bundle, \emph{Compositio Math.} \textbf{101}(1996) 217--224


\bibitem[DPS96]{DPS96}
Demailly, J.-P.; Peternell, T.; Schneider, M. {Holomorphic line
bundles with partially vanishing cohomology.} Proceedings of the
Hirzebruch 65 Conference on Algebraic Geometry (Ramat Gan, 1993),
165--198, Israel Math. Conf. Proc., {\bf 9}, Bar-Ilan Univ 1996.



\bibitem[GHS03]{GHS03} Graber, T.; Harris, J.; Starr, J. Families of rationally
connected varieties. \emph{J. Amer. Math. Soc.} \textbf{16} (2003),
no. 1, 57--67. 

%
%



\bibitem[HW15]{HW15}
 Heier, G; Wong, B.
 On projective {K}\"ahler manifolds of partially positive curvature and rational connectedness.
\href{https://arxiv.org/abs/1509.02149}{arXiv:1509.02149}.



\bibitem[Huy05]{Huy05}Huybrechts, D.\emph{ Complex geometry. An introduction.} Universitext. Springer-Verlag, Berlin,
2005. 




\bibitem[Kod54]{kod54} Kodaira, K. On K\"ahler varieties of restricted type (an intrinsic characterization of algebraic varieties). \emph{Ann. of Math. (2)} \textbf{60}, (1954).
28--48. 


\bibitem[Kol96]{Kol96} Koll\'ar, J. \emph{Rational curves on algebraic
varieties.} Berlin: Springer-Verlag, 1996.


\bibitem[KMM92]{KMM92} Koll\'ar, J.; Miyaoka, Y.; Mori, S.  Rationally connected varieties. \emph{J. Algebraic Geom.} 1 (1992), no. 3,
429--448. 

\bibitem[Lam99]{La99}  Lamari, A.  Courants k\"ahleriens et surfaces compactes. Ann. Inst.
Fourier \textbf{49}(1999), 263--285.

\bibitem[LP17]{LP17} Lazic, V.; Peternell, T. Rationally connected varieties―on a
conjecture of Mumford. \emph{Sci. China Math.} \textbf{60} (2017),
1019--1028.

\bibitem[Liu16]{Liu16} Liu, Gang. Three-circle theorem and dimension estimate for holomorphic functions on K\"ahler manifolds. Duke Math. J. \textbf{165} (2016), no. 15, 2899--2919.


\bibitem[LY17]{LY14} Liu, K.-F.; Yang, X.-K. Ricci curvatures on Hermitian manifolds.{Trans. Amer. Math. Soc.}  \textbf{369} (2017),
5157--5196.

\bibitem[Mat13]{Mat13}
Matsumura, S. {Asymptotic cohomology vanishing and a converse to the
Andreotti-Grauert theorem on surfaces.}  Ann. Inst. Fourier.
Grenoble \textbf{63} (2013) 2199--2221.


\bibitem[Mat18]{Mat18}Matsumura, S. On the image of MRC fibrations of projective manifolds with semi-positive holomorphic sectional
curvature. \href{https://arxiv.org/abs/1801.09081}{arXiv:1801.09081}


\bibitem[Mok88]{Mok88} Mok, Ngaiming. The uniformization theorem for compact K\"ahler manifolds of nonnegative holomorphic bisectional curvature. J. Differential Geom. \textbf{27} (1988), no. 2, 179--214.

\bibitem[Ni98]{Ni98} Ni, Lei. Vanishing theorems on complete K\"ahler manifolds and
their applications. J. Differential Geom. \textbf{50} (1998), no. 1,
89--122.

\bibitem[Ni18]{Ni18}Ni, Lei.  Liouville theorems and a Schwarz Lemma for holomorphic mappings between K\"ahler
manifolds. \href{https://arxiv.org/abs/1807.02674}{arXiv:1807.02674}

\bibitem[NZ18a]{NZ1} Ni, L.; Zheng, F.-Y. Comparison and vanishing theorems for K\"ahler manifolds. \href{https://arxiv.org/abs/1802.08732}{arXiv:1802.08732}

\bibitem[NZ18b]{NZ2} Ni, L.; Zheng, F.-Y. Positivity and Kodaira embedding
theorem. \href{https://arxiv.org/abs/1804.09696}{arXiv:1804.09696}


\bibitem[Ott12]{Ott12} Ottem, J.-C.  Ample subvarieties and $q$-ample divisors.
Adv. Math. \textbf{229} (2012), no. 5, 2868--2887.


\bibitem[Pet06]{Pet06} Peternell, T. Kodaira dimension of subvarieties II.  \emph{Int. J. Math.} \textbf{17} (2006),
619--631. 

\bibitem[Pet12]{Pet12} Peternell, T. Varieties with generically nef tangent bundles. J. Eur. Math. Soc. \textbf{14} (2012), no. 2, 571--603.

%


\bibitem[Tot13]{Tot13}
Totaro, B. {Line bundles with partially vanishing cohomology.}
 \emph{J. Eur. Math. Soc.} \textbf{15} (2013) 731--754.




\bibitem[YZ16]{YZ16}
Yang, B.; Zheng, F.-Y.
 Hirzebruch manifolds and positive holomorphic sectional curvature.
\href{https://arxiv.org/abs/1611.06571v2}{arXiv:1611.06571v2}.






\bibitem[Yang16]{Yang16}Yang, X.-K. Hermitian manifolds with semi-positive
holomorphic sectional curvature.  \emph{Math. Res. Lett.}
\textbf{23} (2016), no.3, 939--952. 



%



\bibitem[Yang17a]{Yang17a}Yang, X.-K. Scalar curvature on compact complex manifolds.
\href{https://arxiv.org/abs/1705.02672}{arXiv:1705.02672}. \emph{To
appear in Trans. Amer. Math. Soc.}



\bibitem[Yang17]{Yang17D}Yang, X.-K. A partial converse to the Andreotti-Grauert
theorem. \href{https://arxiv.org/abs/1707.08006}{arXiv:1707.08006}.
\emph{To appear in Compositio. Math.}

\bibitem[Yang18]{Yang18} Yang, X.-K. RC-positivity, rational connectedness and Yau's conjecture.
{Camb. J. Math.}  \textbf{6} (2018), 183--212.

\bibitem[Yang18a]{Yang18a} Yang, X.-K. RC-positivity, vanishing theorems and rigidity of holomorphic
maps. \href{https://arxiv.org/abs/1807.02601}{arXiv:1807.02601}


\bibitem[Yang18c]{Yang18c}  Yang, X.-K. RC-positivity and
rigidity of harmonic maps into Riemannian manifolds.

\bibitem[Yang2]{Yang2} Yang, X.-K. Rigidity theorems on complete K\"ahler manifolds with
RC-positive curvature.  In preparation.


\bibitem[Yau78]{Yau78} Yau, S.-T.  { On the Ricci curvature of a compact K\"ahler manifold and the complex Monge-Amp\`ere equation, I},
\emph{Comm. Pure Appl. Math.} {\textbf{ 31}} (1978),  339--411.


\bibitem[Yau82]{Yau82} Yau, S.-T. Problem section. In Seminar on
Differential Geometry, \emph{Ann. of Math Stud.} 102, 669-706. 1982.

\end{thebibliography}
\end{document}